\documentclass[a4paper, 10pt]{IEEEtran}

\IEEEoverridecommandlockouts                              
\usepackage{amsmath}
\usepackage{cite}

\usepackage[utf8]{inputenc}
\usepackage{times} 
\usepackage{latexsym,amsfonts,amssymb,amsmath,amscd,mathrsfs}
\usepackage{xargs}   
\usepackage{xspace}
\usepackage{float}
\usepackage{multirow}
\usepackage[tight]{subfigure}
\usepackage{hyperref}
\usepackage{graphicx,tikz,times,pgffor}
\usepackage{pgfplots}
\usetikzlibrary{calc,through,backgrounds,patterns}
\usepackage{tabulary}
\usepackage{tabularx}
\usepackage{changepage}
\usepackage{soul}
\usepackage{todonotes}
\usepackage[ruled,vlined]{algorithm2e}

\newcommand{\setX}{\mathcal{X}}
\newcommand{\X}{\mathbb{X}}

\newcommand{\setE}{\mathcal{E}}

\newcommand{\R}{\mathcal{R}}

\newtheorem {theo}{Theorem}
\newtheorem {propt}{Property}
\newtheorem {defn}{Definition}
\newtheorem {cor}{Corollary}
\newtheorem {lem}{Lemma}
\newtheorem {rem} {Remark}

\newtheorem {pro} {Problem}

\newenvironment{pf}{\noindent \textbf{Proof:}}{}

\title{\LARGE \bf Dynamic characterization of control SIR-type systems and optimal single-interval control}

\author{A. H. Gonz\'alez, A. L. Anderson, A. Ferramosca, E.A. Hernandez-Vargas
\thanks{A.H. González and A. L. Anderson are with the Institute of Technological Development for the Chemical Industry (INTEC), CONICET-Universidad Nacional del Litoral (UNL), Guemes 3450, Santa Fe (3000), Argentina.}%
\thanks{A. Ferramosca is with the Department of Management, Information and Production Engineering, University of Bergamo Via Marconi 5, Dalmine (BG) 24044, Italy.}
\thanks{E.A. Hernandez-Vargas is with the Instituto de Matemáticas, Universidad Nacional Autonoma de Mexico, Boulevard Juriquilla 3001, Santiago de Querétaro, Qro. 76230, Mexico.}%
}

\begin{document}

\maketitle
\thispagestyle{empty}
\pagestyle{empty}

\begin{abstract}
Although modeling studies are focused on the control of SIR-based systems describing epidemic data sets (particularly the COVID-19), few of them present a formal dynamic characterization in terms of equilibrium sets and stability. Such concepts can be crucial to understand not only how the virus spreads in a population, but also how to tailor government interventions such as social distancing, isolation measures, etc. 
The objective of this work is to provide a full dynamic characterization of SIR-type systems under single-interval control actions and, based on it, to find the control action that produces the smallest number of infected individuals at the end of the epidemic that avoids second wave outbreaks. 
Simulations illustrate the benefits of the aforementioned results in terms of the herd immunity threshold.
\end{abstract}

\section{Introduction}

SIR-type models are based on the seminal work of \cite{kermack1927}, which firstly established a compartmental relationship between the main variables of an epidemic: Susceptible, Infected and Removed individuals. After that, several extensions and modifications have been proposed in the literature to account for different epidemic characteristics. Surprisingly, even when significant improvements have been made in the general understanding of the SIR dynamics, the core of almost any epidemic behavior is still qualitatively described by the relatively simple original 3-state model. With the outbreak of the novel COVID-19 (produced by the SARS-CoV-2) at the end of 2019, a plethora of studies have been developed to explain how the virus spread around the world \cite{giordano2020sidarthe}. All of them, with more or less complexity, bases their forecast, analysis and control strategies on SIR-type models. 

In \cite{brauer2012mathematical,sontag2011lecture} a formal analysis was made concerning the general behaviour of SIR-type models and the limit values (for time going to infinity) of the state variables were characterized for several scenarios and initial conditions corresponding to the outbreak of the epidemic. In \cite{harko2014exact} an exact explicit solution for SIR model was given, which allowed the scientific community to go further with a more detailed dynamical analysis. Later on, \cite{franco2020feedback,bertozzi2020challenges} made further analysis concerning the effect of a time varying coefficient with control purposes.   SIR-type models are not only useful to understand epidemic behaviors but also (and more important) to asses the ways one can control them according to \textit{a priori} specified objectives. The qualitative (an previous) evaluation of the different scenarios corresponding to different government measures (social distancing, lockdown, use of face mask, hygiene recommendation, vaccination, etc.) represents a helpful tools for the authorities decision-making in emergency times. At this point, dynamical and control system theory \cite{sontag2013mathematical} becomes crucial to exploit formal mathematical analysis to account for optimal control actions.

The SIDARTHE (Susceptible, Diagnosed, Ailing, Recognized, Threatened, Healed, Extinct) model by \cite{giordano2020sidarthe} was presented to describe the COVID-19 spread in Italy, which provides a detailed description of symptomatic/asymptomatic/isolation/detected state of the individuals, even when the general dynamical behavior is the same of SIR-type models: it shows a critical or threshold value that cannot be surpassed by the final value of the susceptible individuals. 
In \cite{sadeghi2020universal,federico2020taming,morris2021optimal} a rigorous analysis is made of simple social distancing control actions (mainly, the ubiquitous single interval social distancing). By means of different techniques, the authors of these previous works showed how to find the optimal single interval control action, to minimize the infected peak prevalence (maximal fraction of infected individuals). A similar approach concerning single interval interventions is presented in \cite{bliman2021best,di2021optimal}, but to minimize the epidemic final size (or total fraction infected). Most of these proposal consider the rather unrealistic scenario of possible full lockdown, where the transmission rate is zero.
Furthermore, in \cite{kohler2020robust} it is proposed a robust non linear model predictive control (MPC) based on the SIDARTHE model introduced in \cite{giordano2020sidarthe}. The MPC controller manipulates - as usual - the transmission rate (which directly affects the reproduction number), given in this case by several parameters. The objectives are to minimize both, the number of fatalities and the time of isolation, compared to a baseline policy. In \cite{morato2020optimal} an MPC is proposed, based on SIRD (Susceptible-Infected-Recovered-Dead) and SIRASD (Susceptible-Infected-Recovered-Asymptomatic-Symptomatic-Dead) models. A first-order differential equation is used to describe the social distancing, which is a binary signal (on-off, depending on weather the social distancing is or not implemented) affecting the transmission rate. The control objective consists in minimizing both, the current number of infected individuals and the time of isolation. Other similar approaches concerning MPC strategies can be found in \cite{alleman2020covid,peni2020nonlinear,carli2020model}.

The objective of this manuscript is to present a pure dynamical-system perspective to formally analyze SIR-type models and their variants/extensions. The difference between convergence, $\epsilon-\delta$ (Lyapunov) stability and asymptotic stability of the equilibria is established, showing that this kind of analysis is crucial to properly understand the whole system behaviour and the control actions that can be implemented to obtain an optimal closed-loop performance. The key goal of this articles is to solve the following problem:
\begin{pro}
Instead of considering the control objective of minimizing the peak of infected individuals and/or the number of deaths, we consider (as a previous step) to 
maximize the final number susceptible individuals that could avoid second wave outbreaks - which means to minimize the total number of infected individuals in the epidemic.
\end{pro}

The solution of this problem constitutes an important baseline to go further in the design of complex control actions aiming to find the social distancing sequence that also minimizes the infected peak and number of deaths, together with the time the measures last.

\section{Review of SIR Model}\label{sec:revSIR}

This section describes the SIR mathematical model - in its different forms (\cite{kermack1927,brauer2012mathematical,sontag2011lecture}) - from a dynamical system point of view. Consider a population of $N$ individuals in a closed geographical region. The continuous-time dynamical system describing the epidemic is given by the following differential equations
\begin{subequations}\label{eq:SIR}
\begin{align}
    \dot{S}(t) &= - \beta S(t) I(t) \\
    \dot{I}(t) &=   \beta S(t) I(t) - \gamma I(t)  \\
    \dot{C}(t) &= \gamma I(t), 
\end{align}
\end{subequations}
where $S(t)$ is the fractions of individuals who are susceptible to contract the infection at time $t$, $I(t)$ is the fractions of infected/infective individuals (that cause other individuals to become infected) at time $t$, and $C(t)$ is the cumulative fractions of removed individuals (who have already recovered from the disease, or deceased because of it, up to time $t$). Parameters $\beta$ and $\gamma$ stand for the transmission and the recovery/death rates of the disease, respectively. 
\begin{rem}
	Each of the latter compartments ($S$, $I$ and $C$) can be divided into a number of sub-compartments, potentially connected to each others, to have a more detailed description of an epidemic \cite{giordano2020sidarthe}. However, the main dynamic of the original system is, in general, maintained, as detailed in \cite{sadeghi2020universal}.
\end{rem}	

A non-dimensional version of model \eqref{eq:SIR} can be obtained by rescaling the time by $\tau :=t \gamma$ (\cite{sontag2011lecture,bertozzi2020challenges}). This way, model \eqref{eq:SIR} reads:
\begin{subequations}\label{eq:SIRnondim}
\begin{align}
    \dot{S}(\tau) &= - \R S(\tau) I(\tau) \\
    \dot{I}(\tau) &=   \R S(\tau) I(\tau) - I(\tau) \\
    \dot{C}(\tau) &=  I(t),
\end{align}
\end{subequations}
where $\R=\beta/\gamma$ is the so-called \textbf{basic reproduction number}. To emphasize the fact that all the variable are positive, let us define the following constraint set
\begin{eqnarray}
\setX \!\!:=\!\! \{(S,I,C) \in \mathbb R^3\!\!: \! S \in [0,1],\! I \in [0,1],\! C \in [0,1] \}, \nonumber
\end{eqnarray}
in such a way that $(S(\tau),I(\tau),C(\tau)) \in \setX$ for all $\tau \geq 0$. Furthermore, note that $\dot{S}(\tau)+\dot{I}(\tau)+\dot{C}(\tau) = 0$, so $S(\tau)+I(\tau)+R(\tau) = 1$, for $\tau \geq 0$. Particularly, $S(0)+I(0)+R(0) = 1$, where $\tau=0$ is assumed to be the epidemic outbreak time, in such a way that $(S(0),I(0),R(0)):=(1-\epsilon,\epsilon,0)$, with $0 < \epsilon \ll 1$, \textit{i.e.}, the fraction of susceptible individuals is smaller than, but close to $1$; the fraction of infected is close to zero and the fraction of removed are null.

The solution of \eqref{eq:SIRnondim} - which was analytically determined in \cite{harko2014exact}, for $\tau \geq \tau_0 >0$, depends on $\R$ and the initial conditions $(S(\tau_0),I(\tau_0),C(\tau_0))\in \setX$. Since $S(\tau) \geq 0$, $I(\tau) \geq 0$, for $\tau \geq \tau_0 >0$, then $S(\tau)$ is a decreasing function of $\tau$ (by \ref{eq:SIRnondim}.a) and $C(\tau)$ is an increasing function of $\tau$, for all $\tau \geq \tau_0$.
From \eqref{eq:SIRnondim}.b, it follows that if $S(\tau_0)\R \leq 1$, $\dot{I}(\tau) = (\R S(\tau) -1)I(\tau) \leq 0$ at $\tau_0$. Furthermore, given that $S(\tau)$ is decreasing, $I(\tau)$ is also decreasing for all $\tau \geq \tau_0$.
On the other hand, if $S(\tau_0)\R>1$, $I(\tau)$ initially increases, then reaches a global maximum, and finally decreases to zero. In this latter case, the peak of $I$, $\hat I$, is reached at $\hat \tau$, when $\dot I = \R SI - I=0$. This implies that $S=S^*$, where 
\begin{eqnarray}\label{eq:Sstar}
S^*:= \min \{1,1/\R\} 
\end{eqnarray}
is a threshold or critical value, known as "herd immunity". This way, conditions $S(\tau_0)\R>1$ and $S(\tau_0)\R<1$ that determines if $I(\tau)$ increases or decreases at $\tau_0$ can be rewritten as $S(\tau_0)>S^*$ and $S(\tau_0)<S^*$, respectively.\\

For the sake of simplicity, we define $S_\infty:= \lim_{\tau\rightarrow \infty} S(\tau)$, $I_\infty:= \lim_{\tau \rightarrow \infty} I(\tau)$ and $C_\infty:= \lim_{\tau\rightarrow \infty} c(\tau)$, which are values that depend on initial conditions $S(\tau_0),I(\tau_0)$, $C(\tau_0)$, and $\R$. By taking $\tau \rightarrow \infty$ for the solutions proposed in \cite{harko2014exact}, we obtain $I_\infty = 0$, and, so, $C_\infty = 1 - S_\infty$. Furthermore, $S_\infty$ fulfills the condition\footnote{Note that the equation $S_\infty = S(\tau_0) e^{-\R(1-S_\infty)}$ - given in \cite{brauer2012mathematical,franco2020feedback} - is only valid when $C(\tau_0)=0$.}:
\begin{eqnarray}\label{eq:Ssol0}
S(\tau_0)e^{\R C(\tau_0)} e^{-\R(1-S_\infty)} = S_\infty, \nonumber
\end{eqnarray}
which, considering that $S(\tau_0)+I(\tau_0)+C(\tau_0)=1$, can be simplified to
\begin{eqnarray}\label{eq:Ssol}
S(\tau_0)e^{-\R (S(\tau_0)+I(\tau_0))}= S_\infty e^{-\R S_\infty}.
\end{eqnarray}
%



\section{Equilibrium characterization and stability} \label{sec:eq_estabil}

\begin{figure}
	\centering
	\includegraphics[width=0.75\columnwidth]{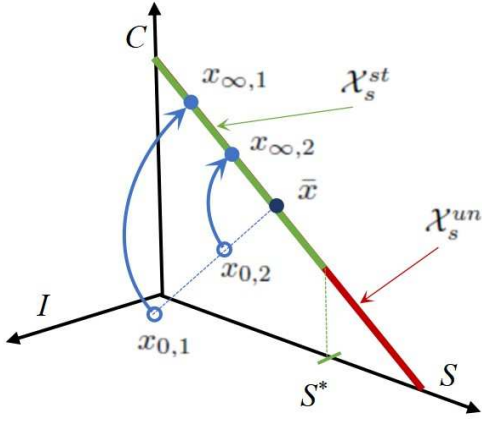}
	\caption{\small{Every point in $\setX_s^{st}$ is $\epsilon-\delta$ stable but not attractive. Initial states $x_0$ starting arbitrarily close to $\bar x \in \setX_s^{st}$ remain (for all $\tau \geq 0$) arbitrarily close to $\bar x$, but do not converge to $\bar x$. As a consequence, set $\setX_s^{st}$ is AS but the points	inside it are not.}}
	\label{fig:stabil_scheme}
	\vspace{-0.3cm}
\end{figure}

The equilibrium of system \eqref{eq:SIRnondim} is obtained by zeroing each of the differential equations. This way the equilibrium set is given by $\setX_s:=\{(\bar S,\bar I,\bar C) \in \mathbb{R}^3: \bar S \in [0,1], \bar I=0, \bar C = 1 - \bar S\}$. For initial conditions $(S(\tau_0),I(\tau_0),C(\tau_0))\in \setX$, the latter set can be refined as 
\begin{eqnarray}
\setX_s \!\!:=\!\!\{(\bar S,\bar I,\bar C) \in \mathbb{R}^3\!\!:\! \bar S \in [0,S(\tau_0)], \bar I\!\!=\!\!0, \bar C \!\!=\!\! 1 \!\!-\!\! \bar S\}, \nonumber
\end{eqnarray}
since $S$ is decreasing.
Next, a key theorem concerning the asymptotic stability of a subset of $\setX_s$ is introduced.
\begin{theo}[Asymptotic Stability]\label{theo:stability}
Consider system \eqref{eq:SIRnondim} constrained by $\setX$. Then, the set
\begin{eqnarray}
\setX_s^{st} \!\!:=\!\! \{(\bar S,\bar I,\bar C) \in \mathbb{R}^3\!\!:\! \bar S \in [0,S^*], \bar I\!\!=\!\!0, \bar C \!\!=\!\! 1 \!\!-\!\! \bar S\}, \nonumber
\end{eqnarray}
where $S^*$ is the herd immunity defined as $S^*:=\min \{1,1/\R\}$, is the unique symptotically stable (AS) of system \eqref{eq:SIRnondim}, with a domain of attraction (DOA) given by $\setX$.  
\end{theo}
\begin{pf}
The proof is divided into two parts. First it is shown that $\setX_s^{st}$ is the smallest attractive equilibrium set in $\setX$. Then, it is shown that $\setX_s^{st}$ is the largest locally $\epsilon-\delta$ stable equilibrium set in $\setX$ which, together with the previous results, implies that $\setX_s^{st}$ is the unique symptotically stable (AS) of system \eqref{eq:SIRnondim}, with a domain of attraction (DOA) given by $\setX$.

\textit{Attractivity}:
%
By manipulating equation \eqref{eq:Ssol}, we have
\begin{eqnarray}\label{eq:Ssol1}
-\R S(\tau_0)e^{-\R (S(\tau_0)+I(\tau_0))}= -\R S_\infty  e^{-\R S_\infty}.
\end{eqnarray}
Denote $y:=-\R S_\infty$ and $z:=-\R S(\tau_0)e^{-\R (S(\tau_0)+I(\tau_0))}$. Then, equation \eqref{eq:Ssol1} can be written as $z=ye^y$. Then $W(z)=y$, where $W(\cdot)$ is the Lambert function. That is, $W(-\R S(\tau_0)e^{-\R (S(\tau_0)+I(\tau_0))}) = -\R S_\infty$, or
\begin{eqnarray}\label{eq:Ssol2}
S_\infty(\R,S(\tau_0),\!I(\tau_0))\!\! :=\!\! -\frac{W(\!-\!\R\! S(\tau_0)\!e^{-\R (S(\tau_0)\!+\!I(\tau_0))})}{\R}.\!\!
\end{eqnarray}
Since $W(z)$ is an increasing function (it goes from $-1$ at $z\!\!=\!\!-1/e$ to $0$ at $z\!\!=\!\!0$), it reaches its minimum at $z\!\!=\!\!-1/e$. $z(S(\tau_0),I(\tau_0))\!\!=\!\!-\R S(\tau_0)e^{-\R (S(\tau_0)\!+\!I(\tau_0))}$ reaches its maximum when $S(\tau_0)\!\!=\!\!S^*$, independently of the values of $\R$ and $I(\tau_0)$ (see Lemmas \ref{lem:Sinf_opt} in the Appendix 2), in which case it is $z(S(\tau_0),I(\tau_0))\!\!=\!\!1/e$. Then, $W(z)$ is bounded from above by $-1$, which means that $S^*$ is an upper bound for $S_\infty$. Therefore, $S_\infty \in [0,S^*]$, which shows the attractivity of $\setX_s^{st}$. 

Figures \ref{fig:SinfFunc} shows a plot of $S_\infty$ as a function of $S(\tau_0)$ and $I(\tau_0)$ for a fixed value of $\R\!>\!1$ (a similar plot can be obtained for $\R\!<\!1$, in which $S_\infty$ reaches its maximum at the vertex of the domain, $S(\tau_0)\!\!=\!\!1$ and $I(\tau_0)\!\!=\!\!0$).
\begin{figure}
	\centering
	\includegraphics[width=0.95\columnwidth]{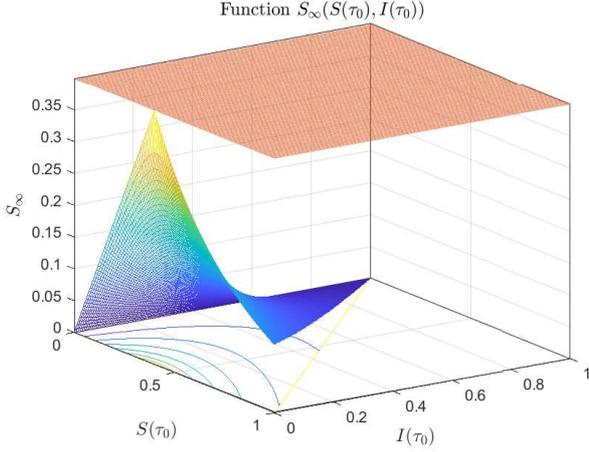}
	\caption{\small{Function $S_\infty(\R,S(\tau_0),C(\tau_0))$, with $\R\!=\!2.5$ is bounded from above by $S^*\!=\!1/\R$ ($S^*\!=\!0.4$, light red plane). Furthermore, $S_\infty$ reaches its maximum, given by $S_\infty^{op} \!=\! S^*$, at $S(\tau_0)\!=\!S^*$, $C(\tau_0) \!=\! 1-S^*$, which implies that $I(\tau_0)\!=\!0$.}}
	\label{fig:SinfFunc}
	\vspace{-0.3cm}
\end{figure}
%
%
%

To show that $\setX_s^{st}$ is the smallest attractive set in $\setX$, consider a state $(\bar S, \bar I, \bar C):= \bar x \in \setX_s^{st}$ and an arbitrary small ball of radius $\epsilon > 0$, w.r.t. $\setX$, around it, $\mathbb{B}_{\epsilon}(\bar x) \in \setX$. Pick two arbitrary initial states $x_{0,1}=(S_{0,1},I_{0,1},C_{0,1})$ and $x_{0,2}=(S_{0,2},I_{0,2},C_{0,2})$ in $\mathbb{B}_{\epsilon}(\bar x)$, such that $S_{0,1} \neq S_{0,2}$ and $C_{0,1} \neq C_{0,2}$. These two states converge, according to equation \eqref{eq:Ssol2}, to $x_{\infty,1}=(S_{\infty,1},0,1-S_{\infty,1})$ and $x_{\infty,2}=(S_{\infty,2},0,1-S_{\infty,2})$, respectively, with $S_{\infty,1},S_{\infty,2} \in [0,S^*]$. Given that function $z(S(\tau_0),I(\tau_0))$ is monotone (injective) in $S(\tau_0)$ and $I(\tau_0)$, and $W(z)$ is monotone (injective) in $z$, then $S_{\infty,1} \neq S_{\infty,2}$. This means that, although both initial states converge to some state in $\setX_s^{st}$, they necessarily converge to different points. Therefore neither single states $\bar x \in \setX_s^{st}$ nor subsets of $\setX_s^{st}$ are attractive in $\setX$, which shows that $\setX_s^{st}$ is the smallest attractive set in $\setX$ (see Figure \ref{fig:stabil_scheme} for a schematic plot of the behavior of single states in $\setX_s^{st}$).

\textit{Local $\epsilon-\delta$ stability}:
Let us consider a particular equilibrium point $\bar x := (\bar S,0,\bar C)$, with $\bar S \in [0,S^*]$ and $\bar C = 1-\bar S$ (i.e., $\bar x \in \setX_s^{st}$).	Then a Lyapunov function candidate is given by (a modified version of the one used in \cite{nangue2019global})
\begin{eqnarray} \label{ec:lya1}
	V(x) := S- \bar S - \bar S \ln(\frac{S}{\bar S}) + I.
\end{eqnarray}
This function is continuous in $\setX$, is positive definite for all non-negative $x \neq \bar x$ and, furthermore, $V(\bar x)=0$. 
Function $V$ evaluated at the solutions of system \eqref{eq:SIRnondim} reads:
\begin{eqnarray} \label{ec:lya2}
\frac{\partial V(x(\tau))}{\partial \tau} \!\!\!\! &=& \frac{\partial V}{\partial x} \dot{x}(\tau) \nonumber\\
&=& \!\!\! \left[\frac{d V}{d S}~\frac{d V}{d I}~ \frac{d V}{d R} \right] \left[
\begin{array}{c}
-\R S(\tau) I(\tau)  \\
\R S(\tau) I(\tau)-I(t)  \\
I(\tau)   
\end{array}\right]\nonumber\\
&=& \!\!\! \left[(1-\frac{\bar S}{S(\tau)})~1~ 0 \right] \left[
\begin{array}{c}
-\R S(\tau) I(\tau)  \\
\R S(\tau) I(\tau)-I(\tau) \\
I(\tau)  
\end{array}\right]\nonumber\\
&=& \!\!\! I(\tau) (\R \bar S  - 1)
\end{eqnarray}
for $x(0)\in \setX$ and $\tau \geq 0$.
Function $\dot{V}(x(\tau))$ depends on $x(\tau)$ only through $I(\tau)$. So, independently of the value of the parameter $\bar S$, $\dot{V}(x(\tau))=0$ for $I(\tau)\equiv 0$. This means that for any single $x(0) \in \setX_s$, $I(0)=0$ and so, $I(\tau)= 0$, for all $\tau\geq 0$. So $\dot{V}(x(\tau))$ is null for any $x(0) \in \setX_s$ (i.e, it is not only null for $x(0)=\bar x$ but for any $x(0) \in \setX_s$).

On the other hand, for $x(0) \notin \setX_s$, function $\dot{V}(x(t))$ is negative, zero or positive, depending on if the parameter $\bar S$ is smaller, equal or greater than  $S^*= \min \{1,1/\R\}$, respectively, and this holds for all $x(0)\in \setX$ and $\tau \geq 0$. So, for any $\bar x \in \setX_s^{st}$, $\dot V(x(\tau))\leq 0$ (particularly, for $\bar x=(\bar S,0,0)=(S^*,0,0)$, $\dot V(x(\tau)) = 0$, for all $x(0) \in \setX$ and $\tau \geq 0$) which means that each $\bar x \in \setX_s^{st}$ is locally $\epsilon-\delta$ stable (see Theorem \ref{theo:lyap} in Appendix 1). Then, if every state in $\setX_s^{st}$ is locally $\epsilon-\delta$ stable, the whole set $\setX_s^{st}$ is locally $\epsilon-\delta$ stable.

Finally, by following similar steps, it can be shown that $\setX_s^{un}$ is not $\epsilon-\delta$ stable, which implies that $\setX_s^{st}$ is also the largest locally $\epsilon-\delta$ stable set in $\setX$, which completes the proof. $\square$\\

\end{pf}
\begin{rem} \label{rem:nonas}
In the latter proof, if we pick a particular $\bar x \in \setX_s^{st}$, then $\dot{V}(x(t))$ is not only null for $x(0)=\bar x$ but for all $x(0) \in \setX_s^{st}$, since in this case, $I(\tau)=0$, for $\tau\geq0$. This means that it is not true that $\dot{V}(x(t)) < 0$ for every $x \neq \bar x$, and this is the reason why single equilibrium points (and subsets of $\setX_s^{st}$) are $\epsilon-\delta$ stable, but not attractive.\\
\end{rem}
Two corollaries of Theorem \ref{theo:stability}, concerning the properties of $\setX_s^{st}$ and the general behavior of $S_\infty$, are presented next.
\begin{cor}\label{cor:Xst}
Consider system \eqref{eq:SIRnondim} with arbitrary initial conditions $(S(\tau_0),I(\tau_0),C(\tau_0))\in \setX$, for some $\tau_0 \geq 0$. Then: 
(i) Set $\setX_s^{st}$ is in general a subset of $\setX_s$ (for $\R<1$, $\setX_s^{st} \equiv \setX_s$), and its size depends on $\R$, but not on the initial conditions.
(ii) Subsets of $\setX_s^{st}$ are $\epsilon-\delta$ stable but not attractive (\textit{i.e.}, even when $\setX_s^{st}$ is AS as a whole, no subset of it is AS). For each single equilibrium state $x_s \in \setX_s^{st}$ there always exists an arbitrary small disturbance (in $I$) such that the state trajectory converges to another state in $\setX_s^{st}$. This is particularly true for the state $(S^*,0,0)$ which belongs to $\setX_s^{st}$.
(iii) If $\R <1$, $S^* = 1$. Then, $\setX_s^{st} \equiv \setX_s$ and the so called healthy equilibrium $x_h:=(\bar S,0,0)$ with $\bar S = 1$, is in $\setX_s^{st}$, and so it is $\epsilon-\delta$ stable, but not attractive (any small value of $I$ will make the system to converge to $(\bar S,0,0)$) with $\bar S < 1$. 
(iv) If $\R > 1$, set $\setX_s$ can be divided into two sets, $\setX_s = \setX_s^{st} \cup \setX_s^{un}$, where 
\begin{eqnarray}
\setX_s^{un}:=\{(\bar S,\bar I,\bar C) \in \mathbb{R}^3: \bar S \in (S^*,1], \bar I=0, \bar C = 1 - \bar S\}, \nonumber
\end{eqnarray}
is \textbf{an unstable equilibrium set} (which contains the healthy equilibrium). 
(v) Given that any compact set including an AS equilibrium set is AS, $\setX_s$ is AS, for any value of $\R$. However, if $\R>1$, it contains an unstable equilibrium set, $\setX_s^{un}$.
\end{cor}

\begin{figure}
	\centering
	\includegraphics[width=0.95\columnwidth]{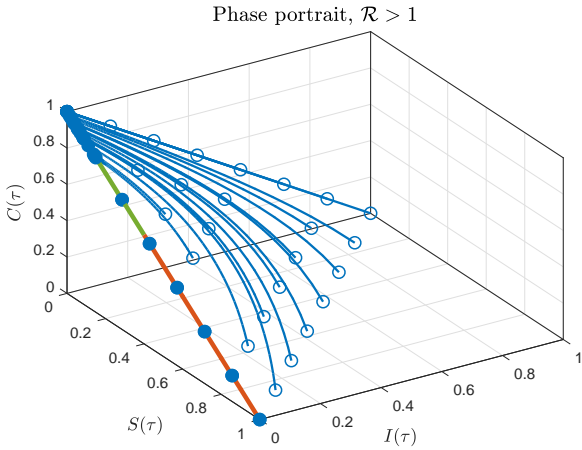}
	\caption{\small{Phase Portrait for system \eqref{eq:SIRnondim} with $\R=2.5$ and starting points such that $S(\tau_0)+I(\tau_0)+C(\tau_0)=1$ (starting points in empty circles, ending points in solid circles). Set $\setX_s^{st}$ is in green, while $\setX_s^{un}$ is in red. As it can be seen, all the trajectories converges to $\setX_s^{st}$. As $\R$ decreases under $1$, $\setX_s^{st}$ is the whole set $\setX_s$.}}
	\label{fig:PhaPorRg1}
	\vspace{-0.3cm}
\end{figure}

Figures \ref{fig:PhaPorRg1} shows a Phase Portrait for system \eqref{eq:SIRnondim}, with $\R\!>\!1$, and initial conditions such that $S(\tau_0)+I(\tau_0)+C(\tau_0)=1$. Similar behavior can be seen if $\R<1$, with $\setX_s^{st} \!\equiv\! \setX_s$.

%
%
%

%
\begin{cor} \label{cor:sinfty}
Consider system \eqref{eq:SIRnondim} with arbitrary initial conditions $(S(\tau_0),I(\tau_0),C(\tau_0)) \in \setX$, for some $\tau_0 \geq 0$. Then:
\begin{enumerate}
\item For any value of $S(\tau_0)>0$ and $I(\tau_0)>0$, $S_\infty(\R,S(\tau_0),I(\tau_0)) \rightarrow 0$ , when $\R \rightarrow \infty$; while $S_\infty(\R,S(\tau_0),I(\tau_0))$ remains close to $S(\tau_0)$ when $\R \rightarrow 0$.
\item For $S(\tau_0) > S^*$ and fixed $I(\tau_0)>0$ and $\R>0$, $S_\infty(\R,S(\tau_0),I(\tau_0))$ decreases with $S(\tau_0)$, and $S_\infty(\R,S(\tau_0),I(\tau_0)) < S^*$. This means that the closer $S(\tau_0)$ is to $S^*$ from above, the closer will be $S_\infty$ to $S^*$ from below.
\item For $S(\tau_0) < S^*$ and fixed $I(\tau_0)>0$ and $\R>0$, $S_\infty(\R,S(\tau_0),I(\tau_0))$ increases with $S(\tau_0)$, and $S_\infty(\R,S(\tau_0),I(\tau_0)) < S^*$. This means that the closer $S(\tau_0)$ is to $S^*$ from below, the further will be $S_\infty$ to $S^*$, from below. 
\item For any fixed $S(\tau_0) \in [0,1]$ and $\R>0$, $S_\infty(\R,S(\tau_0),I(\tau_0))$ decrease with $I(\tau_0)$, and $S_\infty(\R,S(\tau_0),I(\tau_0)) \leq S^*$.
\item For fixed $\R>0$, $S(\tau_0) = S^*$ and $I(\tau_0)=0$, $S_\infty(\R,S(\tau_0),I(\tau_0))$ reaches its maximum over $\{(S,I)\in \mathbb R^2: S\in[0,1],I\in[0,1]\}$, and the maximum value is given by $S^*$ (see Lemma \ref{lem:Sinf_opt} in Appendix 2). 
\end{enumerate}
\end{cor}

The proof of the properties enumerated in Corollary \ref{cor:sinfty} are omitted for brevity. However, Figures \ref{fig:SinfFunc} and \ref{fig:Sinf_S} show how $S_\infty$ behaves for different values of initial conditions.
In the next section the effect of control actions modifying $\R$ for a finite period of time is studied. To simplify the analysis, the following property, concerning quasi steady states, is introduced:
\begin{propt}[Quasi steady state]
Consider system \eqref{eq:SIRnondim} with arbitrary initial conditions $(S(\tau_0),I(\tau_0),C(\tau_0)) \in \setX$, for some $\tau_0 \geq 0$. Then, there is a large enough finite time, $\tau_{qss}$, such that for all $\tau\geq \tau_{qss}$, both $\dot S(\tau)$ and $I(\tau)$ can be approximated by zero, and $S(\tau)$ is approximately equal to $S_\infty(\R,S(\tau_0),I(\tau_0))$. This latter condition is known as Quasi Steady State (QSS) and, opposite to the formal steady state condition, is assumed to be reached in finite time. 
\end{propt}

Time $\tau_{qss}$ depends on the value of $\R$, and can be computed - for practical purposes - as $\tau_{qss}(\R) = 5 \hat \tau(\R)$, where $\hat \tau(\R)$ is peak time of $I(\tau)$ corresponding to a given $\R$.

\section{Control}
Control objectives in epidemic can be defined in several ways. The peak of the infected individual uses to be a critical index to minimize, since it is directly related to the health system capacity. However, other indexes - usually put in a second place - are also important. This is the case of the time the epidemic last, including second (or third) outbreak waves, and the total number of infected individuals at the end of the epidemic (measured by t he epidemic final size $C_\infty=1-S_\infty$). 

Among many others, social distancing is a typical non-pharmaceutical measure that decrease parameter $\beta$ in system \eqref{eq:SIR} or, directly, parameter $\R$ in system \eqref{eq:SIRnondim}. We will consider that $\R$ varies over time but it is piece-wise constant. More precisely, we will assume a single interval social distancing \cite{sadeghi2020universal,federico2020taming,morris2021optimal,bliman2021best,di2021optimal}. At the outbreak of the epidemic ($\tau=0$), it is $(S(0),\!I(0),\!C(0)):=(1\!-\!\epsilon,\epsilon,0)$, with $0 \!<\! \epsilon \ll \! 1$, and $\R(0) \! >\!1$ ($\R(0)$ goes from $2.5$ to $3.5$ for the SARS-CoV-2, all around the world). Then, the single interval social distancing is defined by the following function:
%
%
\begin{eqnarray} \label{eq:sisd}
	\R(\tau) = \left\{ 
	\begin{array}{cc}
		\R(0)  &  \mbox{for}~ \tau \in [0,\tau_s),   \\
		\R_s   & \mbox{for}~ \tau \in [\tau_s,\tau_f],\\
		\R(0)  & \mbox{for}~ \tau \in (\tau_f,\infty),
	\end{array} \right.
\end{eqnarray}
where $\tau_s < \hat \tau(\R(0))$, being $\hat \tau(\R(0))$ the time of the peak of $I(\tau)$ when no social distancing is implemented, $\R_s \in [\R_{min},\R(0)]$, with $\R_{min} \in (0,\R(0))$ (the case $\R_s=0$ is not considered, since it represents a full lockdown which is rather unrealistic) and $\tau_f>\tau_s$, but finite.

The control problem we want to solve reads as follows: for given initial time, $\tau_s < \hat \tau$, find $\R_s$ and $\tau_f$ (finite) to maximize $S_\infty=S_\infty(\R(0),S(\tau_f),I(\tau_f))$. Note that minimize $S_\infty$ corresponds to maximize $C_\infty=1-S_\infty$ (since $I_\infty=0$), which represents the minimal epidemic final size

A critical point concerning social distancing measures (that is usually disregarded) is that they are always temporary control actions ($\tau_f$ finite), not permanent ones (as clearly stated in \cite{sadeghi2020universal} and \cite{kohler2020robust}). It is not possible to maintain efficient social distancing actions for ever (neither for a time long enough to make the virus to disappear) since population fatigue due to psychological or economical problems would systematically relax its effectiveness. 

So, according to the stability results from the previous sections, the following Theorem holds:
\begin{theo}[$S$ steady state upper bound] \label{theo:Suppbound}
	Consider system \eqref{eq:SIRnondim} with initial conditions $(S(0),I(0),C(0))=(1-\epsilon,\epsilon,0)$, $0 < \epsilon \ll 1$, and $\R(0)$ such that $S(0)\!>\!S^*$. Any single interval social distancing as the one defined in \eqref{eq:sisd}, with finite final time, makes the system to converges to an equilibrium state $(S_\infty,0, C_\infty)$ with $S_\infty=S_\infty(\R(0),S(\tau_f),I(\tau_f)) < S^*$, being $S^*<1$ the herd immunity corresponding to no social distancing.
\end{theo}
\begin{pf}
	Under the assumption of a single interval social distancing, Theorem \ref{theo:stability} states that $S_\infty$ is determined by the value of $\R(\tau)$, for  $\tau > \tau_f$, which is $\R(0)$, and by $S(\tau_f),I(\tau_f)$; that is, $S_\infty=S_\infty(\R(0),S(\tau_f),I(\tau_f))$. Furthermore, by Lemma \ref{lem:Sinf_opt}, this function reaches its unique maximum over $\{(S,I)\in \mathbb R^2: S\in[0,1],I\in[0,1]\}$, given by $S^*$, when $I(\tau_f)=0$ and $S(\tau_f)=S^*$. But $I(\tau_f)>0$, given that $\tau_f$ is finite ($I(\tau)$ converges to zero asymptotically), which implies that $S_\infty=S_\infty(\R(0),S(\tau_f),I(\tau_f)) < S^*$, and the proof is complete. $\square$
%
\end{pf}
%

%
\begin{rem}
	Note that even if $R_s$ is time-varying in \eqref{eq:sisd}, Theorem \ref{theo:Suppbound} is still true. This is a simple but strong result concerning any kind of social distancing interrupted at finite time. The minimal possible epidemic final size is completely determined by the epidemic itself (its original $\R(0)$) and, provided that no immunization (by vaccination and/or development of the individuals immune system) is considered, it cannot be modified by non-pharmaceutical measures. The point to be elucidated now is which are the value of $\R_s$ and the finite final time $\tau_f$ that minimize $\|S^*-S_\infty \|$.
\end{rem}

In the search of such a value, the next definition is stated.
\begin{defn}[Goldilocks social distancing]
	The goldilocks social distancing, $\R^{g}=\R^{g}(\tau_s)$, is defined as the one that, if applied at $\tau_s < \hat \tau(\R(0))$, produces $S_\infty(\R^{g},S(\tau_s),I(\tau_s)) = S^*$.
\end{defn}
%

\begin{rem}[$R^g$ computation]
	Given $\tau_s$ and $\R(0)$, $\R^{g} \!\!= \!\!\R^{g}(S(\tau_s),I(\tau_s))\!\!=\!\!\R^{g}(\tau_s)$ can be obtained, numerically, by means of Algorithm \ref{alg:Ropt}. For $0 \!\!< \!\!\tau_s \!\!<\!\! \hat \tau$ and fixed $\R(0)>1$, $\R^{g}(\tau_s)$ is a decreasing function of $\tau_s$.
\end{rem}
	\begin{algorithm}
	\SetAlgoLined
	$\R = \R(0)$\;
	Compute $S(\tau_s)$ and $I(\tau_s)$ by integrating system \eqref{eq:SIRnondim}, with $\R=\R(0)$, from $0$ to $\tau_s$\;
	$S^*=1/\R$\; 
	$S_\infty (\R,S(\tau_s),I(\tau_s)) =$ -lambertw$(0,-\R S(\tau_s) e^{-\R (S(\tau_s)+I(\tau_s))})/\R$\;
	\While{$S_\infty <= S^*$}{
		$\R \leftarrow \R-0.0001$\;
		$S_\infty(\R,S(\tau_s),I(\tau_s)) \leftarrow$ -lambertw$(0,-\R S(\tau_s) e^{-\R (S(\tau_s)+I(\tau_s))})/\R$\;
	}
	$\R^{g}=\R$\;
	\caption{Computation of $\R^{g}(\tau_s)$}
	\label{alg:Ropt}
\end{algorithm}

Clearly, goldilocks social distancing cannot be applied indefinitely, since $\tau_f$ is finite. However, it can be applied up to a time $\tau_f > \tau_s + \tau_{qss}(\R_s) $, in such a way that $S(\tau_f)$ and $I(\tau_f)$ arbitrarily approaches $S^*$ and $0$ (from above), respectively. This way, the following definition can be introduced:
\begin{defn}[Quasi optimal single interval social distancing]\label{teo:cont_act}
	Consider a given starting time, $\tau_s \in (0,\hat \tau(\R(0)))$.
	Then, the \textbf{quasi optimal single interval control action} consists in applying $\R^{g}$, up to a time $\tau_f> \tau_s+\tau_{qss}(\R_s)$, such that the system reaches a QSS (i.e., $S(\tau_f) \approx S^*$ and $I(\tau_f) \approx 0$).
\end{defn}
\begin{rem}
	Clearly, the latter definition refers to a quasi optimal single interval control action, because larger values of $\tau_f$ will produce values of $S(\tau_f)$ and $I(\tau_f)$ closer to $S^*$ and $0$, respectively, so $S_\infty(\R(0),S(\tau_f),I(\tau_f))$ will be closer to $S^*$. 
\end{rem}
\begin{rem}
	If $\R_s$ is allowed to be zero (which, as it was said, is a rather unrealistic case), another quasi optimal single interval control action would be to apply $\R_s=0$ at $\tau_s=\hat{\tau}$ (the time at which $S(\tau)$ reaches $S^*$ and $I(\tau)$ reaches the peak), and keep this control for a large enough time $\tau_f$, such that the system reaches a QSS (i.e., $I(\tau_f) \approx 0$), as proposed in \cite{bliman2021best}.
\end{rem}
\begin{rem}
A question that naturally arises is what happen if, for a given social distancing (quasi-optimal or not), $\tau_f$ cannot be large enough for the system to reach a QSS. Clearly, in this scenario, to implement the goldiclocks social distancing, $\R^{g}$, may be counterproductive, and may produce a value of $S_\infty$ significantly smaller than $S^*$. However, the important point is that all social distancing $\R_s$ implemented for a period of time that not allow the system to reach a QSS before it is interrupted will produce a value of $S_\infty(\R(0),S(\tau_f),I(\tau_f))$ significantly smaller than the one obtained with the quasi optimal single interval control action.  
\end{rem}

The next Theorem, which is one of the main contribution of the work, summarizes the latter results by means of a classification that consider every possible single interval social distancing case.

\begin{figure}
	\centering
	\includegraphics[width=1\columnwidth]{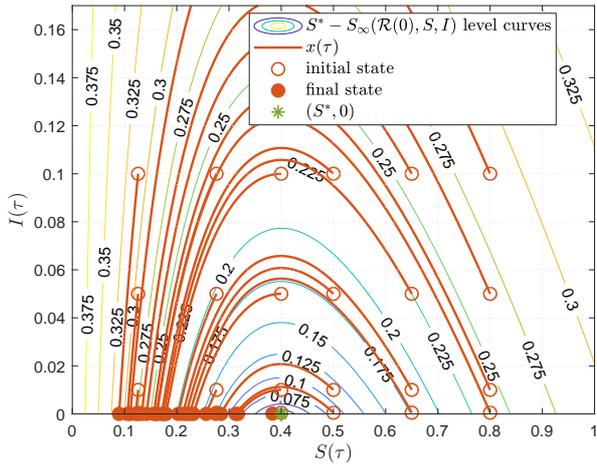}
	\caption{\small{Phase Portrait for system \eqref{eq:SIRnondim} in the $S,I$ plane, with $\R=2.5$ different starting points (red lines), and level curves of function $V(S,I):=S^*-S_\infty(\R,S,I)$. Function $V$ is positive for all $(S,I) \neq (S^*,0)$, is null at $(S^*,0)$ and $\dot V(S,I)=0$ along the solution of system \eqref{eq:SIRnondim} (since $S_\infty(\R,S,I)$ is). So its level sets are arbitrary small invariant sets around $(S^*,0)$. Note that starting states close to $(S^*,0)$ produce time evolution close to $S^*$.}}
	\label{fig:LevCurv}
	\vspace{-0.3cm}
\end{figure}

\begin{theo}[Single interval social distancing scenarios]\label{teo:cont_sce}
	Consider system \eqref{eq:SIRnondim} with initial conditions $(S(0),I(0),C(0))=(1-\epsilon,\epsilon,0)$, $0 < \epsilon \ll 1$, with $\R(0)$ such that $S(0)>S^*$. Consider also single interval social distancing (as the one defined in \eqref{eq:sisd}), with a given starting time $\tau_s \in (0,\hat \tau(\R(0)))$, and a finite final time $\tau_f$. Define soft and strong social distancing depending on if $\R_s>\R^{g}$ or $\R_s<\R^{g}$, respectively. Define also long and short term social distancing depending on if the system reaches or does not reach a QSS at $\tau_f$ (i.e., if $\tau_f>\tau_s + \tau_{qss}(R_s)$ or $\tau_f<\tau_s + \tau_{qss}(R_s)$). Then, 
	the following scenarios can take place:
	\begin{enumerate}
	\item Quasi optimal single interval social distancing: if $\R_s=\R^{g}$, and $\tau_f>\tau_s + \tau_{qss}(R_s)$ (i.e., $(S(\tau_f),I(\tau_f),C(\tau_f))$ reaches a QSS), then $S_\infty(\R(0),S(\tau_f),I(\tau_f)) \approx S^*$. Furthermore, the closer is $S(\tau_f)$ to $S^*$ (or $I(\tau_f)$ to zero), the closer will be $S_\infty(\R(0),S(\tau_f),I(\tau_f))$ to $S^*$.
	\item Soft long term social distancing: if $(S(\tau_f),I(\tau_f),C(\tau_f))$ reaches a QSS, and $S(\tau_f) < S^*$, then $S_\infty(\R(0),S(\tau_f),I(\tau_f)) \approx S(\tau_f) <S^*$; \textit{i.e.}, $S(\tau)$ will remain approximately constant for $\tau \geq \tau_f$. Furthermore, the softer a soft long term social distancing is, the smaller will be $S_\infty(\R(0),S(\tau_f),I(\tau_f))$.
	\item Strong long term social distancing: if $(S(\tau_f),I(\tau_f),C(\tau_f))$ reaches a QSS, and $S(\tau_f) >S^*$, a \textbf{second outbreak} wave will necessarily take place at some time $\hat{\hat{\tau}} > \tau_f$ and, finally, the system will converge to an $S_\infty(\R(0),S(\tau_f),I(\tau_f))< S^*$. Furthermore, the stronger a strong long term social distancing is, the larger will be the second wave and the smaller will be $S_\infty(\R(0),S(\tau_f),I(\tau_f))$. 
	\item Short term social distancing: if $(S(\tau_f),I(\tau_f),C(\tau_f))$ does not reach a QSS (\textit{i.e}, if $I(\tau_f) \not\approx 0$), then soft, strong and goldilocks social distancing will necessarily produce values of $S_\infty(\R(0),S(\tau_f),I(\tau_f))$ significantly smaller than the one obtained by quasi optimal single interval social distancing. In general, larger values of $I(\tau_f)$ will produce smaller values of $S_\infty(\R(0),S(\tau_f),I(\tau_f))$. 
	This case includes the particular case where the social distancing is interrupted at the very moment at which $S(\tau_f)=S^*$, but with $I(\tau_f) \not\approx 0$. This means that the herd immunity needs to be reached a steady state, not as a transitory one.
	%
	\end{enumerate}
\end{theo}
\begin{pf}
	The proof follows from the stability results shown in Section \ref{sec:eq_estabil}, and equation \eqref{eq:Ssol2} defining $S_\infty$:
	\begin{enumerate}
		\item Given that $\R_s=\R^g$ is implemented for $\tau\in[\tau_s,\tau_f]$, $\tau_f$ is finite but greater than $\tau_s + \tau_{qss}(R_s)$ and $S_\infty (\R^g,S(\tau_s),I(\tau_s)) = S^*$, then $S(\tau_f)$ approaches $S^*$ and $I(\tau_f)$ approaches zero, from above, as $\tau_f$ increases. This means that at $\tau_f$, when social distancing is interrupted, $(S(\tau_f),I(\tau_f),C(\tau_f))$ is close to the unstable equilibrium set $\setX_s^{un}$. Then, by Corollary \ref{cor:sinfty}.(2), function $S_\infty(\R(0),S(\tau_f),I(\tau_f))$ is such that the closer $(S(\tau_f),I(\tau_f),C(\tau_f))$ is to the equilibrium point $(S^*,0,1-S^*)$, with $S(\tau_f)>S^*$, the closer will be $S_\infty(\R(0),S(\tau_f),I(\tau_f))$ to $S^*$, with $S_\infty<S^*$ (see the 'pine' shape of $S_\infty$ around $S^*$, for $I\approx 0$, in Figure \ref{fig:Sinf_S})\footnote{Indeed, by the $\epsilon-\delta$ stability of the equilibrium state $(S^*,0,1-S^*)$, for each (arbitrary small) $\epsilon>0$, it there exists $\delta>0$, such that, if the system starts in a ball of radius $\delta$ centered at $(S^*,0,1-S^*)$, it will keeps indeterminately in the ball of radius $\epsilon$ centered at $(S^*,0,1-S^*)$. Furthermore, it is possible to define invariant sets around $(S^*,0,1-S^*)$ by considering the level sets of the Lyapunov function \eqref{ec:lya1}, with $\bar S=S^*$, or even the level sets of function $V(S,I,C):=S^*- S_\infty(\R,S,I)$, with a fixed $\R>0$. This way, once the system enters any arbitrary small level set of the latter functions, it cannot leaves the set anymore. See, Figure \ref{fig:LevCurv}}.
		\item Given that $(S(\tau_f),I(\tau_f),C(\tau_f))$ approaches a steady state with $S(\tau_f) < S(\tau_f)$, then $(S(\tau_f),I(\tau_f),C(\tau_f))$ is close to the stable equilibrium set $\setX_s^{st}$, when the social distancing is interrupted. Then, the system will converge to an equilibrium with $S_\infty(\R(0),S(\tau_f),I(\tau_f))$ close to $S(\tau_f)$. Softer social distancing produces smaller values of $S(\tau_f)$ and, by Corollary \ref{cor:sinfty}.(3), smaller values of $S(\tau_f)$ produce smaller values of $S_\infty(\R(0),S(\tau_f),I(\tau_f))$. 
		\item Given that $(S(\tau_f),I(\tau_f),C(\tau_f))$ approaches a steady state with $S(\tau_f) > S^*$, then $(S(\tau_f),I(\tau_f),C(\tau_f))$ is close to the unstable equilibrium set, $\setX^{un}$, when the social distancing is interrupted. Then, the system will converges to an equilibrium in the stable equilibrium set, $\setX^{st}$, with $S_\infty(\R(0),S(\tau_f),I(\tau_f)) <S^*$. Stronger social distancing produces greater values of $S(\tau_f)$ and, by Corollary \ref{cor:sinfty}.(2), values of $S(\tau_f)$ farther from $S^*$, from above, produce values of $S_\infty(\R(0),S(\tau_f),I(\tau_f))$ farther from $S^*$, from below. When $S(\tau_f)$ is significantly greater than $S^*$, no matter how large is $\tau_f$ and how small is $I(\tau_f)$\footnote{Note that as long as $\tau_f$ is finite, $(S(\tau_f),I(\tau_f),C(\tau_f))$ cannot reach $\setX_s^{un}$, and so $I(\tau_f)$, even when arbitrary small, is greater than zero. So, once the social distancing is interrupted, the system evolves to an equilibrium in $\setX_s^{un}$. But even if we assume that for a large $\tau_f$, $(S(\tau_f),I(\tau_f),C(\tau_f))$ reaches $\setX_s^{un}$ (the last infected individual is recovered), an infected individual entering the population from outside, will destabilize the system in the same way.}, the system will evolve to an equilibrium in $\setX_s^{st}$, with $S_\infty(\R(0),S(\tau_f),I(\tau_f))$ significantly smaller than $S^*$. Furthermore, to go from $S(\tau_f)$ to $S_\infty(\R(0),S(\tau_f),I(\tau_f))$, for $\tau>\tau_f$, the system significantly increase $I(\tau)$, and this effect is known as a second outbreak wave.
		\item Given that $(S(\tau_f),I(\tau_f),C(\tau_f))$ is a transitory state, then it does not approach any equilibrium. This means that $I(\tau_f)$ is significantly greater than $0$, and according to Lemma \ref{lem:Sinf_opt}, in Appendix 2, the maximum of $S_\infty(\R(0),S(\tau_f),I(\tau_f))$ over $\setE(\delta)\!\!:=\!\!\{(S(\tau_f),I(\tau_f)) \!\!\in\!\! \mathbb R^2\!\!:\!\! S(\tau_f)\!\!\in\!\! [0,1],~I(\tau_f) \!\!\in\!\! [\delta,\!1]\}$ is given by $- W(- \R S^* e^{-\R (S^*+\delta)})/\R$, which a decreasing function of $\delta$, and reaches $S^*$ only when $\delta=0$ (see Figure \ref{fig:Sinf_S}). Then, independently of the value of $S(\tau_f)$, $S_\infty(\R(0),S(\tau_f),I(\tau_f))$ will be significantly smaller than the one obtained with optimal single interval social distancing, in which $\delta \approx 0$. $\square$
	\end{enumerate}
\end{pf}
%


\subsection{Some remarks}
%
(i) The quasi optimal single interval social distancing does not account for the peak of the infected individuals, but only for the maximal final value of $S$ that avoids further outbreaks. Indeed, if a piece-wise constant $\R(\tau)$ is applied for $\tau \in [\tau_s,\tau_f]$ (instead of a single value) the peak of $I$ can be arbitrary selected to be under any threshold (defined, for instance, by ICU system).\\
(ii) To apply soft control actions (even no-social-distancing at all), expecting it would evolve alone to the herd immunity value, is not an option. $S_\infty$ is in general significantly smaller than $S^*$ (specially for the reported values of $\R$ for the COVID-19).\\
(iii) To apply intermediate social distancing up to the time the system reaches the herd immunity, and then interrupt it, is not an option. This strategy implies to reduce the reproduction number not so drastically, in such a way that $S$ underpasses $S^*$ at some time. However, if the social distancing is interrupted at a transient state (assuming $\R_s>0$), the initial conditions for the next time period are such that $S_\infty$ will be significantly smaller than $S^*$.\\
(iv) To apply hard social distancing measures, for an undefined large period of time, expecting the epidemic will die out alone is not an option. Hard social distancing produces an $S(\tau_f)$ larger than $S^*$ (no matter if $\R_s>1$ or $\R_s<1$), but this $S_f$ is an artificial steady state value, since once the social distancing is dropped or reduced, a second infection wave will necessarily occurs at some future time.\\
(v) The results presented in the previous sections are quite universal for all kind of SIR-type models describing the COVID-19 epidemic. Indeed, as long as a reproduction number is well defined and three groups of compartments are stated - as in the SIDHARTE model presented in \cite{giordano2020sidarthe} - the stability and optimality results are still valid.
\section{Simulation Results}
To illustrate the results of the previous sections, some simulations are performed next. Consider system \eqref{eq:SIRnondim} with $\R=2.5$ (which corresponds to the herd immunity $S^*=0.4$), and initial conditions$(S(0),I(0),C(0))=(0.9950,0.0050,0.0000)$. Consider also a single interval social distancing starting at $\tau_s = 2$ (in adimensional units). Then, different values of $\R_s$ (the social distancing reproduction number), and $\tau_f$ (final single interval social distancing time) will be considered to demonstrate the different cases described in Theorem \ref{teo:cont_sce}.

\subsection{Long term social distancing}
 
\begin{figure}
	\centering
	\includegraphics[width=0.95 \columnwidth]{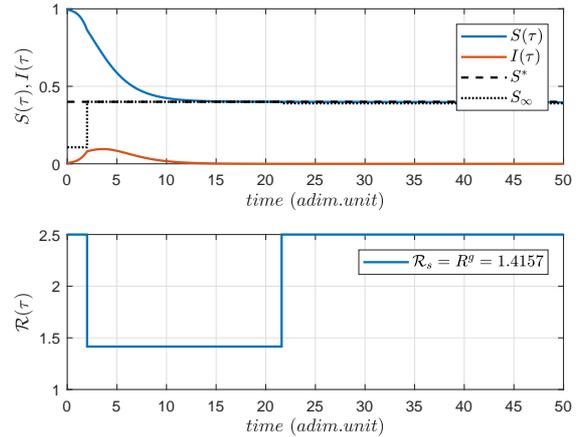}
	\caption{\small{System evolution corresponding to the quasi optimal single interval social distancing $\R_s=\R^g=1.4157$, interrupted after a quasi steady state is reached. $S_\infty(\R(0),S(\tau_f),I(\tau_f)) = 0.3942$, which as predicted, is a value close to $S^*=0.4$.}}
	\label{fig:QOp}
	\vspace{-0.4cm}
\end{figure}

\begin{figure}
	\centering
	\includegraphics[width=0.95 \columnwidth]{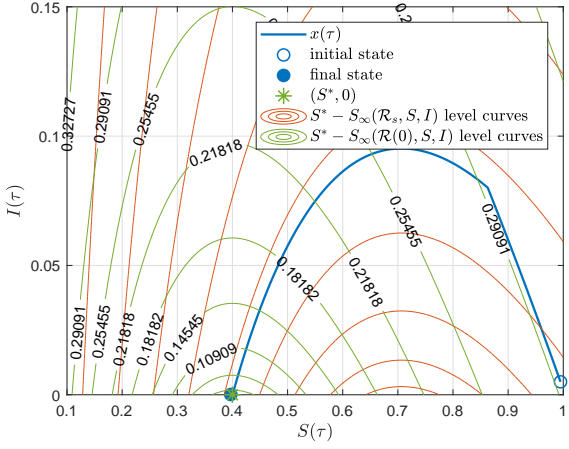}
	\caption{\small{Phase Portrait for system \eqref{eq:SIRnondim} in the $S,I$ plane, with $\R_s=\R^g=1.4157$, and level curves of functions $V(S,I):=S^*-S_\infty(\R_s,S,I)$ (in red) and $V(S,I):=S^*-S_\infty(\R(0),S,I)$ (in green).}}
	\label{fig:QOppp}
	\vspace{-0.4cm}
\end{figure}

Consider the case of single interval social distancing during a large enough time for the system to reach a QSS before it is interrupted. To ensure such a QSS condition, $\tau_f$ is computed as $6$ times the infected peack time corresponding to quasi optimal single interval social distancing, i.e., $\tau_f = 6 \hat \tau(\R_s)=6(3.6)=21.6$ (in adimensional units). First, the quasi optimal single interval social distancing, which consists in applying $\R_s=\R^g=1.4157$ from $\tau_s$ to $\tau_f$ is simulated. The results are shown in Figure \ref{fig:QOp}, where it can be seen that the social distancing is interrupted when $S(\tau_f)$ is close to $S^*$ (from above) and, so, $S_\infty(\R(0),S(\tau_f),I(\tau_f))$ ends up at a value close to $S^*$, from below. Furthermore, in Figure \ref{fig:QOppp}, the phase portrait of the simulation is shown, in the plane $S,I$, together with the level curves of the Lyapunov functions $V(S,I):=S^*-S_\infty(\R,S,I)$, with $\R$ corresponding to both, the open loop value $\R(0)$ and the social distancing value, $\R_s$. As it can be seen, the state goes along the open loop Lyapunov level curves (in red) up to time $\tau_s$, then switches to a trajectory along the social distancing Lyapunov level curves (in green) up to $\tau_f$, to finally return to a trajectory along the open loop Lyaponv level curves. This later trajectory cannot be seen in this case, since at $\tau_f$ the system is very close to $\epsilon-\delta$ stable equilibrium $(S^*,0)$.

\begin{figure}
	\centering
	\includegraphics[width=0.95 \columnwidth]{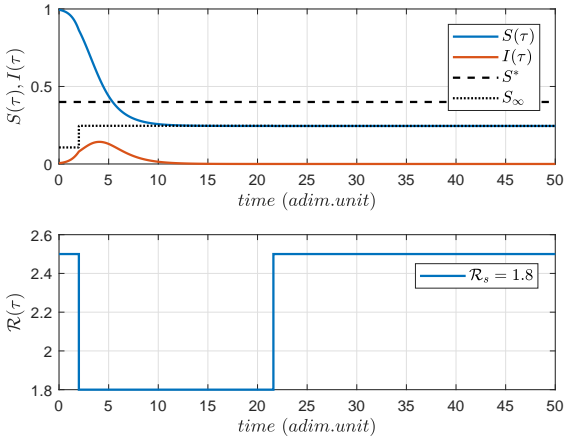}
	\caption{\small{System evolution corresponding to a social distancing $\R_s=1.8$, interrupted after a quasi steady state is reached. $S_\infty(\R(0),S(\tau_f),I(\tau_f)) = 0.2453$, which as predicted, is a value significantly smaller than $S^*=0.4$.}}
	\label{fig:SO}
	\vspace{-0.4cm}
\end{figure}

\begin{figure}
	\centering
	\includegraphics[width=0.95 \columnwidth]{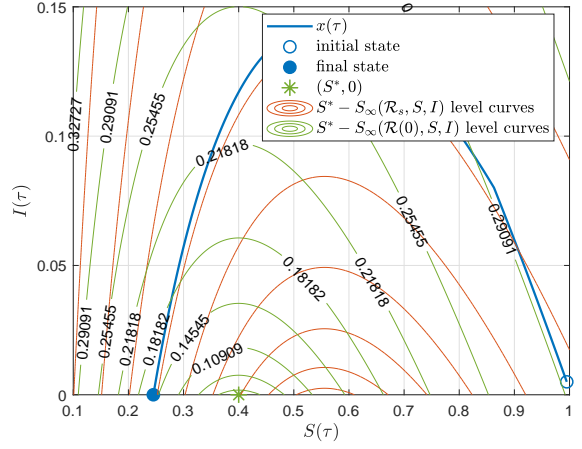}
	\caption{\small{Phase Portrait for system \eqref{eq:SIRnondim} in the $S,I$ plane, with $\R_s=1.8$, and level curves of functions $V(S,I):=S^*-S_\infty(\R_s,S,I)$ (in red) and $V(S,I):=S^*-S_\infty(\R(0),S,I)$ (in green).}}
	\label{fig:SOpp}
	\vspace{-0.4cm}
\end{figure}

Next, a soft social distancing is considered, by applying $\R_s=1.8$ from $\tau_s$ to $\tau_f$. The results are shown in Figure \ref{fig:SO}, where it can be seen that the social distancing is interrupted when $S(\tau_f)$ is close to its QSS value, $S_\infty(\R_s,S(\tau_s),I(\tau_s))$, which is smaller than $S^*$. So, $S_\infty(\R(0),S(\tau_f),I(\tau_f))$ ends up at a value close to $S_\infty(\R_s,S(\tau_s),I(\tau_s))$, which as it was said, is significantly smaller than $S^*$. Furthermore, in Figure \ref{fig:SOpp}, the phase portrait of the simulation is shown, in the plane $S,I$, together with the level curves of the Lyapunov functions $V(S,I):=S^*-S_\infty(\R,S,I)$, with $\R$ corresponding to both, the open loop value $\R(0)$ and the social distancing value, $\R_s$. The state goes along the open loop Lyapunov level curves (in red) up to time $\tau_s$, then switches to a trajectory along the social distancing Lyapunov level curves (in green) up to $\tau_f$, to finally return to a trajectory along the open loop Lyaponv level curves. This later trajectory cannot be seen in this case, since at $\tau_f$ the system is very close to $\epsilon-\delta$ stable equilibrium $(S_\infty(\R(0),S(\tau_f),I(\tau_f)),0)$.

\begin{figure}
	\centering
	\includegraphics[width=0.95 \columnwidth]{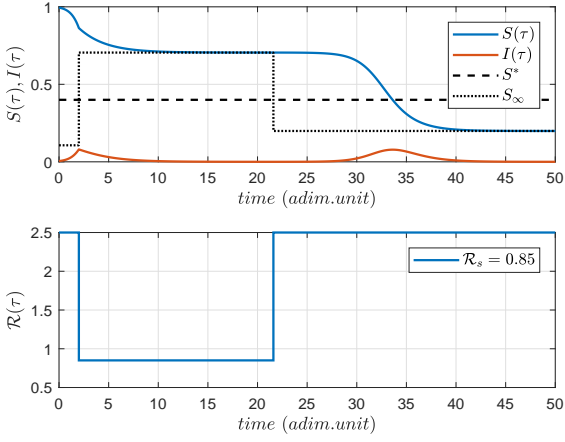}
	\caption{\small{System evolution corresponding to a social distancing $\R_s=0.85$, interrupted after a quasi steady state is reached. $S_\infty(\R(0),S(\tau_f),I(\tau_f)) = 0.1989$, which as predicted, is a value significantly smaller than $S^*=0.4$.}}
	\label{fig:SW}
	\vspace{-0.4cm}
\end{figure}

Finally, a strong social distancing is considered, by applying $\R_s=0.85$ from $\tau_s$ to $\tau_f$. The results are shown in Figure \ref{fig:SW}, where it can be seen that the social distancing is interrupted when $S(\tau_f)$ is close to its QSS value, $S_\infty(\R_s,S(\tau_s),I(\tau_s))$, which is significantly greater than $S^*$. So, $S_\infty(\R(0),S(\tau_f),I(\tau_f))$ ends up at a value significantly smaller than $S^*$, after producing a significantly large infected second wave. Furthermore, in Figure \ref{fig:QOppp}, the phase portrait of the simulation is shown, in the plane $S,I$, together with the level curves of the Lyapunov functions $V(S,I):=S^*-S_\infty(\R,S,I)$, with $\R$ corresponding to both, the open loop value $\R(0)$ and the social distancing value, $\R_s$. The state goes along the open loop Lyapunov level curves (in red) up to time $\tau_s$, then switches to a trajectory along the social distancing Lyapunov level curves (in green) up to $\tau_f$, to finally return to a trajectory along the open loop Lyaponv level curves. This later trajectory is significant (with a large excursion of $I(\tau)$), since it goes from $S_\infty(\R_s,S(\tau_s),I(\tau_s)) = 0.7$, to $S_\infty(\R(0),S(\tau_f),I(\tau_f))=0.1989$. Note that the state $(S(\tau_f),I(\tau_f))$ is close to an $\epsilon-\delta$ unstable equilibrium, and no matter how large is $\tau_f$, the system will experience a second wave of the same magnitude.
\begin{figure}
	\centering
	\includegraphics[width=0.95 \columnwidth]{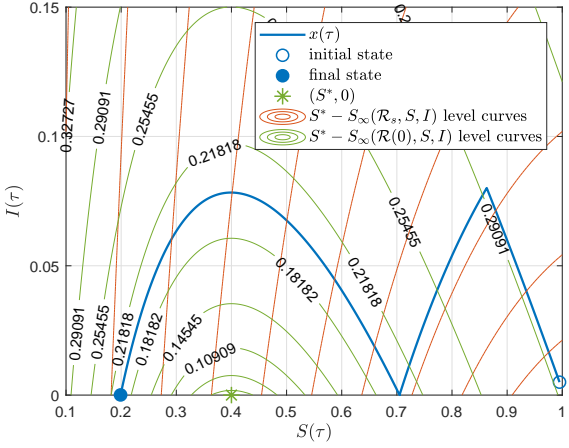}
	\caption{\small{Phase Portrait for system \eqref{eq:SIRnondim} in the $S,I$ plane, with $\R_s=0.85$, and level curves of functions $V(S,I):=S^*-S_\infty(\R_s,S,I)$ (in red) and $V(S,I):=S^*-S_\infty(\R(0),S,I)$ (in green).}}
	\label{fig:SWpp}
	\vspace{-0.4cm}
\end{figure}

\subsection{Short term social distancing}

\begin{figure}
	\centering
	\includegraphics[width=0.95 \columnwidth]{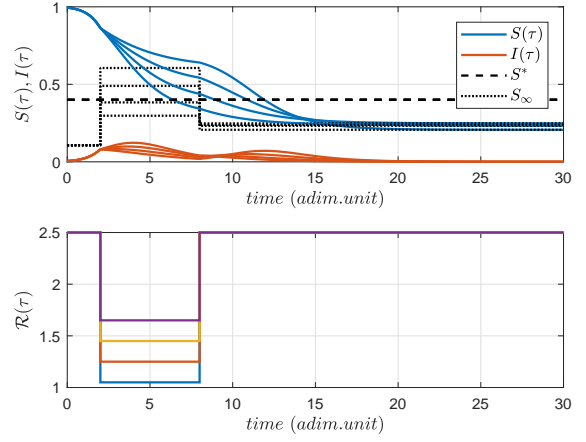}
	\caption{\small{System evolution corresponding to a social distancing $\R_s$ ranging from $0.85$ to $1.8$, interrupted before a quasi steady state is reached. $S_\infty(\R(0),S(\tau_f),I(\tau_f))$ assume the values $0.2066$, $0.2322$, $0.2480$, $0.2384$, whose maximum ($0.2480$) is significantly smaller than the one obtained with the quasi optimal single interval social distancing $0.3942$.}}
	\label{fig:TR}
	\vspace{-0.4cm}
\end{figure}

\begin{figure}
	\centering
	\includegraphics[width=0.95 \columnwidth]{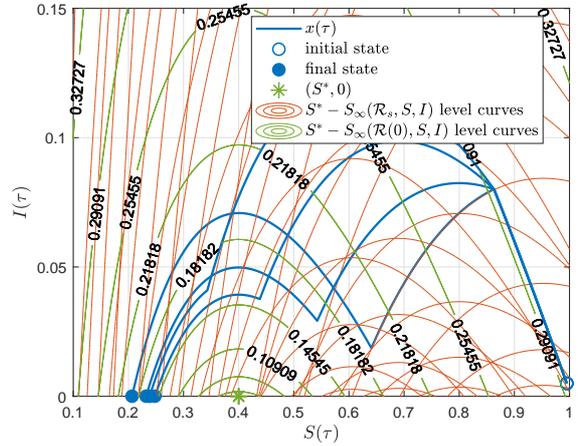}
	\caption{\small{Phase Portrait for system \eqref{eq:SIRnondim} in the $S,I$ plane, with several values of $\R_s$, and level curves of functions $V(S,I):=S^*-S_\infty(\R_s,S,I)$ (in red) and $V(S,I):=S^*-S_\infty(\R(0),S,I)$ (in green).}}
	\label{fig:TRpp}
	\vspace{-0.4cm}
\end{figure}

Now, the case in which the social distancing is interrupted at a transient regime (i.e., when $I(\tau_f) \not\approx 0$) is simulated. To do that, $\tau_s=2$, as before, but $\tau_f$ is reduced to $\tau_f=8$ (in adimensional units). Furthermore, several values of $\R_s$ (all producing $I(\tau_f) \not\approx 0$, for $\tau_f=8$) are considered, to show that $S_\infty(\R(0),S(\tau_f),I(\tau_f))$ is always significantly smaller that value obtained with the quasi optimal single interval social distancing. Figure \ref{fig:TR} shows the simulation for $\R_s$ ranging from $1.05$ to $1.8$ (close to the soft and strong values of the latter subsection). In Figure \ref{fig:TRpp}, it is shown the phase portrait of the simulation in the plane $S,I$, together with the level curves of the Lyapunov functions $V(S,I):=S^*-S_\infty(\R,S,I)$, with $\R$ corresponding to both, the open loop value $\R(0)$ and the social distancing values, $\R_s$. As it can be seen, the system does not approach a steady state at time $\tau_f$, so the Lyapunov level curves along with the system evolves after $\tau_f$ (in green), steers the state to values far from $(S^*,0)$.
Note that even when the social distancing is interrupted when $S(\tau_f)\approx S^*$, a small value of $S_\infty(\R(0),S(\tau_f),I(\tau_f))$ is obtained, since $I(\tau_f) \not\approx 0$.

\section{CONCLUSIONS AND FUTURE WORKS}

In this work, the stability and general long term behavior of SIR-type models were discuses in the context of controlling the COVID-19 epidemic. A quasi optimal control action - consisting in the single period social distancing producing the greater final fraction of susceptible individuals - is found. Several suboptimal scenarios are also analyzed.

\section*{Appendix 1. Stability theory}\label{sec:app1}

All the following definitions are referred to system 
\begin{eqnarray}\label{eq:difeqini}
	\dot x (t) = f(x(t)),~~x(0)=x_0,
\end{eqnarray}
where $x$ is the system state constrained to be in $\X \subseteq \R^n$, $f$ is a Lipschitz continuous nonlinear function,
and $\phi(t;x)$ is the solution for time $t$ and initial condition $x$.
\begin{defn}[Equilibrium set]
	Consider system \ref{eq:difeqini} constrained by $\X$. The set $\setX_s \subset \X$ is an equilibrium set if each point $x \in \setX_s$ is such that $f(x) =  0$ (this implying that $\phi(t;x)  =  x$ for all $t \geq 0$).
\end{defn}
\begin{defn}[Attractivity of an equilibrium set]\label{def:attrac_set}
	Consider system \ref{eq:difeqini} constrained by $\X$ and a set $\setX \subseteq \X$. A closed equilibrium set $\setX_s \subset \setX$ is attractive in $\setX$	if $\lim_{t \rightarrow \infty} \|\phi(t;x)\|_{\setX_s} =0$ for all $x \in \setX$.
\end{defn}

A closed subset of an attractive set (for instance, a single equilibrium point) is not necessarily attractive. On the other hand, any set containing an attractive set is attractive, so the significant attractivity concept in a constrained system is given by the smallest one.

\begin{defn}[Local $\epsilon - \delta$ stability of an equilibrium set]\label{def:eps_del_stab}
	Consider system \ref{eq:difeqini} constrained by $\X$. A closed equilibrium set $\setX_s \subset \X$ is $\epsilon-\delta$ locally stable
	if for all $\epsilon >0$ it there exists $\delta>0$ such that in a given boundary of $\setX_s$, $\|x\|_{\setX_s} <\delta$, it
	follows that $ \|\phi(t;x)\|_{\setX_s} < \epsilon$, for all $t \geq 0$.
\end{defn}

Unlike attractive sets, a set containing a locally $\epsilon-\delta$ stable equilibrium set is not necessarily locally $\epsilon-\delta$ stable.

\begin{defn}[Asymptotic stability (AS) of an equilibrium set]\label{def:AS}
	Consider system \ref{eq:difeqini} constrained by $\X$ and a set $\setX \subseteq \X$. A closed equilibrium set $\setX_s \in \X$ is asymptotically stable (AS) in $\setX$ if it is $\epsilon-\delta$ locally stable and attractive in $\setX$.
\end{defn}

Next, the theorem of Lyapunov, which refers to single equilibrium points and provides sufficient conditions for both, local $\epsilon-\delta$  stability and assymptotic stability, is introduced.
\begin{theo}[Lyapunov theorem \textcolor{blue}{\cite{khalil2002nonlinear}}]\label{theo:lyap}
	Consider system \ref{eq:difeqini} constrained by $\X$ and an equilibrium state $x_s \in \X$. Let consider a 
	function $V(x): \R^n \rightarrow\R$ such that $V(x)>0$ for $x \neq x_s$, $V(x_s)=0$ and $\dot{V}(x(t)) \leq 0$, denoted as Lyapunov function.
	Then, the existence of such a function in a boundary of $x_s$ implies that $x_s \in \setX_s$ is locally $\epsilon-\delta$ stable. If in addition
	$\dot{V}(x(t)) < 0$ for all $x\in \setX \subseteq \X$ such that $x \neq x_s$, and $\dot{V}(x_s) = 0$, then $x_s$ is asymptotically stable in $\setX$.
\end{theo}
%

\section*{Appendix 2. Maximum of $S_\infty(\R,S,I)$}\label{sec:app2}

\begin{lem}[Maximum of the function $S_\infty(\R,S,I)$]\label{lem:Sinf_opt}
	Consider the function $S_\infty(\R,S,I) \!\!=\!\! -\frac{W(-\R S e^{-\R (S+I)})}{\R}$, with $\R$ fixed. Then, $S_\infty^{op}(\delta)\!\!:=\!\!\max\limits_{S,I} \{S_\infty(\R,S,I) \!\!:\!\! (S,I) \in \setE(\delta) \} \!\!=\!\! - W(- \R S^* e^{-\R (S^*+\delta)})/\R$, where $\setE(\delta)\!\!:=\!\!\{(S,I) \!\!\in\!\! \mathbb R^2\!\!:\!\! S\!\!\in\!\! [0,1],~I \!\!\in\!\! [\delta,\!1]\}$ is a set of initial conditions with $I\!\!\geq\!\!\delta$, for some fixed $\delta \!\!\in\!\! [0,1]$, and $S^*\!\!:=\!\!1/\R$. Furthermore, $(S^{op}(\delta),\!I^{op}(\delta))\!\!:=\!\!arg\max\limits_{S,I} \{S_\infty(\R,S,I) \!\!:\!\! (S,I) \!\!\in\!\! \setE \} \!\!=\!\! (S^*,\delta)$. Finally, the maximal value of $S_\infty^{op}(\delta)$ over $\delta \!\!\in \!\![0,1]$ is given by $S_\infty^{op}\!\!:=\!\! \max\limits_{\delta} \{S_\infty^{op}(\delta) \!\!:\!\! \delta \in [0,1] \}\!\!=\!\!S_\infty^{op}(0)\!\!=\!\!S^*$.
\end{lem}
\begin{pf}
	According to \eqref{eq:Ssol2}, $S_\infty$ is given by 
	\begin{eqnarray} \label{eq:Sinf1}
	S_\infty(\R,S,I):= - W(-f(\R,S,I))/\R, 
	\end{eqnarray}
	with $f(\R,S,I):=\R S e^{-\R (S+I)}$. Given that $-W(-x)$ is an increasing (injective) function of $x \in [0,1/e]$ and $\R$ is fixed, then $S_\infty(\R,S,I)$ achieves its maximum over $\setE(\delta)$ at the same values of $S$ and $I$ as $f(\R,S,I)$ (next it is shown that $f(\R,S,I) \in [0,1/e]$ for all $(S,I) \in \setE(\delta)$) and $\delta \in [0,1]$. Then, we focus the attention in finding the maximum (and the maximizing variables) of $f(\R,S,I)$. Let denotes the maximum of $f$ as $f^{op}(\delta):=\max\limits_{S,I} \{f(S,I) : (S,I) \in \setE(\delta) \}$, while the maximizing variables are $S^{op}(\delta)$ and $I^{op}(\delta)$.
	
	Given that the maximum of $f$ occurs at the minimal values of $I$, let us consider, for simplicity, that $g(S,I) = I-\delta$, in such a way that we want to solve $(S^{op}(\delta),I^{op}(\delta))=arg\max \{f(\R,S,I) : g(S,I) \leq 0\}$ (we ignore the conditions $0 \leq S \leq 1$ and $I \leq 1$, but it is easy to see the no maximum is achieved at the boundaries of these constraints). Then $\bigtriangledown f =[\frac{\partial f}{\partial S},~\frac{\partial f}{\partial I}] = [\R e^{-\R (S+I)}(1-\R S),~ \R^2 S e^{-\R (S+I)}]$ and $\bigtriangledown g =[\frac{\partial g}{\partial S},~\frac{\partial g}{\partial I}] = [0,~1]$. Optimality conditions can be written as $\bigtriangledown f = \lambda \bigtriangledown g$, where $\lambda \in \mathbb R_{\geq0}$ is a Lagrange multiplier. Then, $[\R e^{-\R (S^{op}(\delta)+I^{op}(\delta))}(1-\R S^{op}(\delta)),~ \R^2 S^{op}(\delta) e^{-\R (S^{op}(\delta)+I^{op}(\delta))}]=[0,~\lambda]$, which implies that
	$\R e^{-\R (S^{op}(\delta)+I^{op}(\delta))}(1-\R S^{op}(\delta))=0$ and $\R^2 S^{op}(\delta) e^{-\R (S^{op}(\delta)+I^{op}(\delta))} = \lambda$. Since $\R>0$, the first equality implies that $1-\R S^{op}(\delta) =0$, or $S^{op}(\delta) = \min\{1/\R\}=S^*$ (since $S^{op}(\delta)$ must be in $[0,1]$). This way, the second equality reads $\R^2 S^{*} e^{-\R (S^{*}+I^{op}(\delta))} = \lambda$, which is true for any value of $I^{op}(\delta)\in [\delta,1]$ and $\lambda>0$. As we know that larger values of $f$ are achieved for smaller values of $I$, then, $I^{op}(\delta)=\delta$.

	The maximum value of $S_\infty$ is then given by $S_\infty^{op}(\delta) = S_\infty(S^{op}(\delta),I^{op}(\delta))$, which reads
	\begin{eqnarray} \label{eq:Sinfmax}
	S_\infty^{op}(\delta)&=& - W(-\R S^{op}(\delta) e^{-\R (S^{op}(\delta)+I^{op}(\delta))})/\R \nonumber\\
						 &=& - W(-\R S^* e^{-\R (S^*+\delta)})/\R.
	\end{eqnarray}
	If $\R\!\!\geq\!\!1$, then $\R S^* \!\!=\!\!1$, and so $S_\infty^{op}(\delta)\!\!=\!\! - W(- e^{-\R (S^*+\delta)})/\R$. On the other hand, if $\R\!\!<\!\!1$, then $\R S^* \!\!=\!\!\R$ and $S_\infty^{op}(\delta)\!\!=\!\! - W(- \R e^{-\R (S^*+\delta)})/\R$.
	
	In any case, $S_\infty^{op}(\delta)$ is a decreasing function of $\delta \!\!\in\!\! [0,1]$, which means that $S_\infty^{op}\!\!:=\!\! \max\limits_{\delta} \{S_\infty^{op}(\delta) \!\!:\!\! \delta \!\!\in\!\! [0,1] \}\!\!=\!\!S_\infty^{op}(0)$, which is given by $S_\infty^{op}(0)= - W(- e^{-1})/\R = 1/\R = S^*$,
	%
	%
	if $\R\geq1$, and $S_\infty^{op}(0) =- W(-\R e^{-\R})/\R = \R/\R =1 =S^*$,
		%
	%
	if $\R<1$. This concludes the proof. $\square$
\end{pf}

\begin{figure}
	\centering
	\includegraphics[width=0.9\columnwidth]{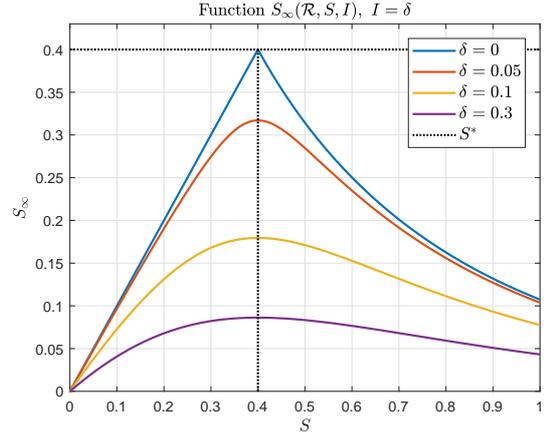}
	\caption{\small{Function $S_\infty(\R,S,I)$, with fixed reproduction number $\R=2.5$, for different values of $\delta$, when $S \in [0,1]$ and $I = \delta$. As it can be seen, the maximum of $S_\infty(\R,S,\delta)$, for any fixed $\delta$, is achieved at $S=S^*$. Particularly, the maximum of $S_\infty(\R,S,\delta)$ over $\delta \in [0,1]$ is achieved at $\delta=0$ and it is given by $S^*$.}}
	\label{fig:Sinf_S}
	\vspace{-0.4cm}
\end{figure}

Figure \ref{fig:Sinf_S} shows function $S_\infty(S,I)$ for different values of $\delta$, when $I =\delta$ and the reproduction number is fixed at $\R=2.5$. Clearly, the maximum is achieved at $S=S^*$, when $\delta=0$. 
\begin{rem}
	In Lemma \ref{lem:Sinf_opt}, parameter $\delta \in [0,1]$ represents the minimal initial value for the fraction of infected individuals. This way, $\delta>0$ means that the initial conditions are not an equilibrium (an equilibrium implies $I=0$). 
\end{rem}

%




\bibliographystyle{IEEEtran}
\bibliography{BSIR_Stabil}

\begin{thebibliography}{10}
\providecommand{\url}[1]{#1}
\csname url@samestyle\endcsname
\providecommand{\newblock}{\relax}
\providecommand{\bibinfo}[2]{#2}
\providecommand{\BIBentrySTDinterwordspacing}{\spaceskip=0pt\relax}
\providecommand{\BIBentryALTinterwordstretchfactor}{4}
\providecommand{\BIBentryALTinterwordspacing}{\spaceskip=\fontdimen2\font plus
\BIBentryALTinterwordstretchfactor\fontdimen3\font minus
  \fontdimen4\font\relax}
\providecommand{\BIBforeignlanguage}[2]{{%
\expandafter\ifx\csname l@#1\endcsname\relax
\typeout{** WARNING: IEEEtran.bst: No hyphenation pattern has been}%
\typeout{** loaded for the language `#1'. Using the pattern for}%
\typeout{** the default language instead.}%
\else
\language=\csname l@#1\endcsname
\fi
#2}}
\providecommand{\BIBdecl}{\relax}
\BIBdecl

\bibitem{kermack1927}
W.~O. Kermack and A.~G. McKendrick, ``A contribution to the mathematical theory
  of epidemics,'' \emph{Proceedings of the royal society of london. Series A,
  Containing papers of a mathematical and physical character}, vol. 115, no.
  772, pp. 700--721, 1927.

\bibitem{giordano2020sidarthe}
G.~Giordano, F.~Blanchini, R.~Bruno, P.~Colaneri, A.~Di~Filippo, A.~Di~Matteo,
  M.~Colaneri \emph{et~al.}, ``A sidarthe model of covid-19 epidemic in
  italy,'' \emph{arXiv preprint arXiv:2003.09861}, 2020.

\bibitem{brauer2012mathematical}
F.~Brauer and C.~Castillo-Chavez, \emph{Mathematical models for communicable
  diseases}.\hskip 1em plus 0.5em minus 0.4em\relax SIAM, 2012.

\bibitem{sontag2011lecture}
E.~D. Sontag, ``Lecture notes on mathematical systems biology,'' 2011.

\bibitem{harko2014exact}
T.~Harko, F.~S. Lobo, and M.~Mak, ``Exact analytical solutions of the
  susceptible-infected-recovered (sir) epidemic model and of the sir model with
  equal death and birth rates,'' \emph{Applied Mathematics and Computation},
  vol. 236, pp. 184--194, 2014.

\bibitem{franco2020feedback}
E.~Franco, ``A feedback sir (fsir) model highlights advantages and limitations
  of infection-based social distancing,'' \emph{arXiv preprint
  arXiv:2004.13216}, 2020.

\bibitem{bertozzi2020challenges}
A.~L. Bertozzi, E.~Franco, G.~Mohler, M.~B. Short, and D.~Sledge, ``The
  challenges of modeling and forecasting the spread of covid-19,''
  \emph{Proceedings of the National Academy of Sciences}, vol. 117, no.~29, pp.
  16\,732--16\,738, 2020.

\bibitem{sontag2013mathematical}
E.~D. Sontag, \emph{Mathematical control theory: deterministic finite
  dimensional systems}.\hskip 1em plus 0.5em minus 0.4em\relax Springer Science
  \& Business Media, 2013, vol.~6.

\bibitem{sadeghi2020universal}
M.~Sadeghi, J.~Greene, and E.~Sontag, ``Universal features of epidemic models
  under social distancing guidelines,'' \emph{bioRxiv}, 2020.

\bibitem{federico2020taming}
S.~Federico and G.~Ferrari, ``Taming the spread of an epidemic by lockdown
  policies,'' \emph{Journal of Mathematical Economics}, p. 102453, 2020.

\bibitem{morris2021optimal}
D.~H. Morris, F.~W. Rossine, J.~B. Plotkin, and S.~A. Levin, ``Optimal,
  near-optimal, and robust epidemic control,'' \emph{Communications Physics},
  vol.~4, no.~1, pp. 1--8, 2021.

\bibitem{bliman2021best}
P.-A. Bliman and M.~Duprez, ``How best can finite-time social distancing reduce
  epidemic final size?'' \emph{Journal of theoretical biology}, vol. 511, p.
  110557, 2021.

\bibitem{di2021optimal}
F.~Di~Lauro, I.~Z. Kiss, and J.~C. Miller, ``Optimal timing of one-shot
  interventions for epidemic control,'' \emph{PLOS Computational Biology},
  vol.~17, no.~3, p. e1008763, 2021.

\bibitem{kohler2020robust}
J.~K{\"o}hler, L.~Schwenkel, A.~Koch, J.~Berberich, P.~Pauli, and
  F.~Allg{\"o}wer, ``Robust and optimal predictive control of the covid-19
  outbreak,'' \emph{arXiv preprint arXiv:2005.03580}, 2020.

\bibitem{morato2020optimal}
M.~M. Morato, S.~B. Bastos, D.~O. Cajueiro, and J.~E. Normey-Rico, ``An optimal
  predictive control strategy for covid-19 (sars-cov-2) social distancing
  policies in brazil,'' \emph{arXiv preprint arXiv:2005.10797}, 2020.

\bibitem{alleman2020covid}
T.~Alleman, E.~Torfs, and I.~Nopens, ``Covid-19: from model prediction to model
  predictive control,'' \emph{https://biomath. ugent.
  be/sites/default/files/2020-04/Alleman\_etal\_v2. pdf, accessed April},
  vol.~30, p. 2020, 2020.

\bibitem{peni2020nonlinear}
T.~P{\'e}ni, B.~Csutak, G.~Szederk{\'e}nyi, and G.~R{\"o}st, ``Nonlinear model
  predictive control with logic constraints for covid-19 management,''
  \emph{Nonlinear Dynamics}, vol. 102, no.~4, pp. 1965--1986, 2020.

\bibitem{carli2020model}
R.~Carli, G.~Cavone, N.~Epicoco, P.~Scarabaggio, and M.~Dotoli, ``Model
  predictive control to mitigate the covid-19 outbreak in a multi-region
  scenario,'' \emph{Annual Reviews in Control}, 2020.

\bibitem{nangue2019global}
A.~Nangue, ``Global stability analysis of the original cellular model of
  hepatitis c virus infection under therapy,'' \emph{American Journal of
  Mathematical and Computer Modelling}, vol.~4, no.~3, pp. 58--65, 2019.

\bibitem{khalil2002nonlinear}
H.~K. Khalil and J.~W. Grizzle, \emph{Nonlinear systems}.\hskip 1em plus 0.5em
  minus 0.4em\relax Prentice hall Upper Saddle River, NJ, 2002, vol.~3.

\end{thebibliography}

\end{document}